\documentclass[11pt]{article}
\usepackage{amsmath,amssymb,geometry,verbatim,algorithm2e,enumerate,float,cite,setspace}
\newtheorem{theorem}{Theorem}

\newtheorem{lemma}[theorem]{Lemma}

\usepackage{setspace}
\usepackage{lineno}
\begin{document}
\onehalfspace
%\linenumbers

\title{Two Greedy Consequences for Maximum Induced Matchings}

\author{Dieter Rautenbach\\[3mm]
\normalsize Institut f\"{u}r Optimierung und Operations Research,
Universit\"{a}t Ulm, Ulm, Germany\\
\normalsize \texttt{dieter.rautenbach@uni-ulm.de}
}

\date{}

\maketitle

\begin{abstract}
We prove that, for every integer $d$ with $d\geq 3$, 
there is an approximation algorithm 
for the maximum induced matching problem restricted to $\{ C_3,C_5\}$-free $d$-regular graphs 
with performance ratio $0.708\bar{3}d+0.425$,
which answers a question posed by
Dabrowski et al. (Theor. Comput. Sci. 478 (2013) 33-40).
Furthermore, we show that every graph with $m$ edges 
that is $k$-degenerate and of maximum degree at most $d$ with $k<d$,
has an induced matching with at least $m/((3k-1)d-k(k+1)+1)$ edges.

\bigskip

\end{abstract}

\pagebreak

\section{Introduction}

A set $M$ of edges of a graph $G$ is an {\it induced matching} of $G$
if the set of vertices of $G$ that are incident with the edges in $M$
induces a $1$-regular subgraph of $G$,
or, equivalently, 
if $M$ is an independent set of the square of the line graph of $G$.
The {\it induced matching number} $\nu_2(G)$ of $G$ 
is the maximum cardinality of an induced matching of $G$.
Induced matchings were introduced by Stockmeyer and Vazirani \cite{stva}
as a variant of classical matchings \cite{lopl}.
While classical matchings are structurally and algorithmically well understood \cite{lopl},
it is hard to find a maximum induced matching \cite{stva,ca}
and efficient algorithms are only known for special graph classes \cite{brmo,ca2,dadelo}.
The problem to determine a maximum induced matching in a given graph,
called {\sc Maximum Induced Matching} for short,
is even APX-complete for bipartite $d$-regular graphs for every $d\geq 3$ \cite{dumazi,dadelo}.

On the positive side, a natural greedy strategy applied to a $d$-regular graph $G$, 
which mimics the well-known greedy algorithm for the maximum independent set problem 
applied to the square of the line graph of $G$,
produces an induced matching with at least $\frac{m(G)}{2d^2-2d+1}$ edges.
Since every induced matching of a $d$-regular graph contains at most $\frac{m(G)}{2d-1}$ edges,
this already yields an approximation algorithm
for {\sc Maximum Induced Matching} in $d$-regular graphs with performance ratio $d-\frac{1}{2}+\frac{1}{4d-2}$
as observed by Zito \cite{zi}.
This was improved slightly by Duckworth et al. \cite{dumazi} 
who describe an approximation algorithm with asymptotic performance ratio $d-1$.
The best known approximation algorithm for {\sc Maximum Induced Matching} 
restricted to $d$-regular graphs 
for general $d$
is due to Gotthilf and Lewenstein \cite{gole}
who elegantly combine a greedy strategy with a local search algorithm
to obtain a performance ratio of $0.75d+0.15$.
In \cite{jorasa} Joos et al. describe a linear time algorithm that finds an induced matching with at least $\frac{m(G)}{9}$ edges for a given $3$-regular graph $G$,
which yields an approximation algorithm for cubic graphs with performance ratio $\frac{9}{5}$.

At the end of \cite{dadelo} Dabrowski et al. propose to study 
approximation algorithms for regular bipartite graphs,
and to determine whether the above performance ratios 
can be improved in the bipartite case.
As our main result we show that this is indeed possible.
\begin{theorem}\label{theorem1}
For every integer $d$ with $d\geq 3$, 
there is an approximation algorithm 
for {\sc Maximum Induced Matching} restricted to $\{ C_3,C_5\}$-free $d$-regular graphs 
with performance ratio $\frac{17}{24}d+\frac{17d}{48d-24}\leq 0.708\bar{3}d+0.425$.
\end{theorem}
Our proof of Theorem \ref{theorem1}
builds on the approach of Gotthilf and Lewenstein \cite{gole},
and all proofs are postponed to Section \ref{section2}.

Our second result,
which also relies on a greedy strategy, 
is a lower bound on the induced matching number of degenerate graphs.
This result is related to recent bounds on the {\it strong chromatic index} 
$\chi_s'(G)$ of a graph $G$ \cite{fj}, 
which is defined as the minimum number of induced matchings 
into which the edge set of $G$ can be partitioned.
The most prominent conjecture concerning this notion 
was made by Erd\H{o}s and Ne\v{s}et\v{r}il in 1985
and states that the strong chromatic index of a graph $G$ of maximum degree at most $d$ 
is at most $\frac{5}{4}d^2$.
The most significant progress towards this conjecture 
is due to Molloy and Reed \cite{more} who proved
$\chi_s'(G)\leq 1.998d^2$
provided that $d$ is sufficiently large.
Again a natural greedy edge coloring implies $\chi_s'(G)\leq 2d^2-2d+1$.

Recall that a graph $G$ is {\it $k$-degenerate} for some integer $k$,
if every non-empty subgraph of $G$ has a vertex of degree at most $k$.
Recently, Chang and Narayanan \cite{cn} studied 
the strong chromatic index of $2$-degenerate graphs
and inspired the following results 
about the strong chromatic index of a $k$-degenerate graph $G$ of maximum degree at most $d$
with $k\leq d$:
\begin{eqnarray}\label{edegen}
\chi_s'(G) \leq
\left\{
\begin{array}{ll}
10d-10 & \mbox{, if $k=2$ \cite{cn}}\\
8d-4 & \mbox{, if $k=2$ \cite{ly}}\\
(4k-1)d-k(2k+1)+1 & \mbox{, \cite{dgs}}\\
(4k-2)d-k(2k-1)+1 & \mbox{, \cite{y}}\\
(4k-2)d-2k^2+1 & \mbox{, \cite{w}}.
\end{array}
\right.
\end{eqnarray}
The proofs of (\ref{edegen}) all rely in some way on greedy colorings 
and (\ref{edegen}) immediately implies that $\nu_s(G)\geq\frac{m(G)}{\chi_s'(G)}\geq \frac{m(G)}{4kd+O(k+d)}$
for the considered graphs.
We show that the factor $4$ can be reduced to $3$.
\begin{theorem}\label{theorem2}
If $G$ is a $k$-degenerate graph of maximum degree at most $d$ with $k<d$,
then $\nu_s(G)\geq \frac{m(G)}{(3k-1)d-k(k+1)+1}.$
\end{theorem}
Before we proceed to the proofs of Theorems \ref{theorem1} and \ref{theorem2}
we collect some notation and terminology.

We consider finite, simple, and undirected graphs, and use standard terminology and notation.
For a graph $G$, we denote the vertex set, edge set, order, and size by
$V(G)$, $E(G)$, $n(G)$, and $m(G)$, respectively.
If $G$ has no cycle of length $3$ or $5$, then $G$ is $\{ C_3,C_5\}$-free.
Let $L(G)$ denote the line graph of $G$
and let $G^2$ denote the square of $G$.
For an edge $e$ of $G$, let 
$C_G(e)=\{ e\}\cup N_{L(G)^2}(e)=\{ f\in E(G):{\rm dist}_{L(G)}(e,f)\leq 2\}$
and let
$c_G(e)=|C_G(e)|$.
Note that a set of edges of $G$ is an induced matching
if and only if it does not contain two distinct edges $e$ and $f$
with $f\in C_G(e)$, or, equivalently, $e\in C_G(f)$.
In a maximal induced matching $M$ of a graph $G$,
for every edge $f$ in $E(G)\setminus M$,
there is some edge $e$ in $M$
with $f\in C_G(e)$.
For some edges $f$, the choice of $e$ might be unique,
which motivates the following definition.
For a set $M$ of edges of $G$, let
$PC_G(M,e)=C_G(e)\setminus \bigcup_{f\in M\setminus \{ e\}}C_G(f)$
and let
$pc_G(M,e)=|PC_G(M,e)|$.
For two disjoint sets $X$ and $Y$ of vertices of $G$,
let $E_G(X,Y)$ be the set of edges $uv$ of $G$ with $u\in X$ and $v\in Y$,
and let $m_G(X,Y)=|E_G(X,Y)|$.
For a set $E$ of edges of $G$, 
let $G-E=(V(G),E(G)\setminus E)$.
For a positive integer $k$, let $[k]=\{ 1,2,\ldots,k\}$.

\section{Proofs}\label{section2}

The greedy strategy for maximum induced matching relies on the following lemma.

\begin{lemma}\label{lemma1}
Let $G$ be a graph. 
If $G_0,\ldots,G_k$ are such that 
\begin{itemize}
\item $G_0=G$, and 
\item for $i\in [k]$, 
there is an edge $e_i$ of $G_{i-1}$ such that 
$G_i=G_{i-1}-C_{G_{i-1}}(e_i)$.
\end{itemize}
then,
\begin{enumerate}[(i)]
\item If $e$ and $f$ are edges of $G_i$ for some $i\in [k]$, then 
$f\in C_{G_i}(e)$ if and only if $f\in C_G(e)$.
\item The set $\{ e_1,\ldots,e_k\}$ is an induced matching.
\end{enumerate}
\end{lemma}
{\it Proof:} 
(i) Let $e=uv$ and $f=xy$ be edges of $G_i$ for some $i\in [k]$.
Since $G_i$ is a subgraph of $G$,
it follows that
$f\in C_{G_i}(e)$ immediately implies $f\in C_G(e)$.
Hence, for a contradiction, we may assume that 
$f\in C_G(e)\setminus C_{G_i}(e)$.
Since $f\not\in C_{G_i}(e)$, the two edges $e$ and $f$ do not share a vertex.
Since $f\in C_G(e)$, we may assume that the graph $G$ contains the edge $ux$.
Since $f\not\in C_{G_i}(e)$, the edge $ux$ does not belong to $G_i$, 
which implies $ux\in C_{G_{j-1}}(e_j)$ for some $j\leq i$.
Note that $e$ and $f$ are edges of $G_{j-1}$.
If $ux$ is incident with $e_j$,
then one of $e$ and $f$ is incident with $e_j$,
and belongs to $C_{G_{j-1}}(e_j)$.
If $ux$ is not incident with $e_j$,
then, by symmetry, 
we may assume that $G_{j-1}$ contains an edge between $u$ and a vertex incident with $e_j$,
which implies that $e$ belongs to $C_{G_{j-1}}(e_j)$.
In both cases we obtain the contradiction 
that one of the two edges $e$ and $f$ does not belong to $G_j$
and hence also not to $G_i$.

(ii) If $e_j\in C_G(e_i)$ for some $i,j\in [k]$ with $i<k$, then, by (i),
we have $e_j\in C_{G_{i-1}}(e_i)$,
which implies the contradiction that $e_j$ does not belong to $G_i$
and hence also not to $G_{j-1}$. 
$\Box$

\bigskip

\noindent Algorithm \ref{alg1}, called {\sc Greedy$(f)$}, 
corresponds to the greedy algorithm used by Gotthilf and Lewenstein.

\bigskip

\begin{algorithm}[H]
{\sc Greedy$(f)$}\\
\KwIn{A $d$-regular graph $G$.}
\KwOut{A pair $(M,G')$ such that $M$ is an induced matching of $G$ and $G'$ is a subgraph of $G$.}
\BlankLine
$G_0\leftarrow G$\;
$G'\leftarrow G_0$\;
$i\leftarrow 1$\;
$M\leftarrow \emptyset$\;
\While{$\min\{ c_{G_{i-1}}(e):e\in E(G_{i-1})\}\leq f$}
{
Choose an edge $e_i$ of $G_{i-1}$ with $c_{G_{i-1}}(e_i)\leq f$\;
$M\leftarrow M\cup \{ e_i\}$\;
$G_i\leftarrow G_{i-1}-C_{G_{i-1}}(e_i)$\;
$G'\leftarrow G_i$\;
$i\leftarrow i+1$\;
}
\Return $(M,G')$\;
\
\caption{The greedy algorithm of Gotthilf and Lewenstein depending on $f$.}\label{alg1}
\end{algorithm}

\bigskip

\noindent Obviously, {\sc Greedy$(f)$} can be performed in polynomial time, 
and, by Lemma \ref{lemma1}, it works correctly.
The intuitive idea behind {\sc Greedy$(f)$} 
is to restrict the greedy choices to edges that are not too expensive
in the sense that their inclusion in the induced matching 
does not eliminate too many of the remaining edges.
In fact, {\sc Greedy$(f)$} adds a fraction of at least $\frac{1}{f}$ 
of the edges in $E(G)\setminus E(G')$ to $M$,
which is better than the trivial fraction $\frac{1}{2d^2-2d+1}$ for $f<2d^2-2d+1$.
The drawback of {\sc Greedy$(f)$} is that it might not consume all edges of $G$,
that is, $E(G)\setminus E(G')$ might be small compared to $E(G)$.

By Lemma \ref{lemma1}, the union of $M$ with any induced matching of $G'$
is an induced matching of $G$.
Therefore Gotthilf and Lewenstein combine {\sc Greedy$(f)$} applied to $G$
with Algorithm \ref{alg2}, called {\sc Local Search}, applied to $G'$.
{\sc Local Search} starts with an empty matching $M'$ 
and performs the following two simple augmentation operations as long as possible:
\begin{itemize}
\item Add an edge from $E(G')\setminus M'$ to $M'$ if this results in an induced matching.
\item Replace one edge in $M'$ with two edges from $E(G')\setminus M'$ if this results in an induced matching.
\end{itemize}
Clearly, {\sc Local Seach} can be performed in polynomial time,
and, as observed above, 
the union of $M$ with its output $M'$ is an induced matching of $G$.
The crucial observation is that the graph $G'$ produced by {\sc Greedy$(f)$} 
has an additional structural property.
As a subgraph of $G$, it is trivially a graph of maximum degree at most $d$
but additionally each of its edges $e$ satisfies $c_{G'}(e)>f$,
which allows the improved analysis of {\sc Local Seach}
in the following two Lemmas \ref{lemma2} and \ref{lemma3}. 

\bigskip

\begin{algorithm}[H]
{\sc Local Search}\\
\KwIn{A graph $G$.}
\KwOut{An induced matching $M$ of $G$.}
\BlankLine
$M\leftarrow \emptyset$\;
\Repeat{$M$ does not increase during one iteration}
%\While{the set $M\cup \{ e\}$ is an induced matching of $G$ for some edge $e\in E(G)\setminus M$, or the set 
%$(M\setminus \{ e\})\cup \{ e',e''\}$ is an induced matching of $G$ 
%for some three distinct edges $e\in M$ and $e',e''\in E(G)\setminus M$}
{
\If{$M\cup \{ e\}$ is an induced matching of $G$ for some edge $e\in E(G)\setminus M$}
{
$M\leftarrow M\cup \{ e\}$\;
}
\If{$(M\setminus \{ e\})\cup \{ e',e''\}$ is an induced matching of $G$ 
for some three distinct edges $e\in M$ and $e',e''\in E(G)\setminus M$}
{
$M\leftarrow (M\setminus \{ e\})\cup \{ e',e''\}$\;
}
}
\Return $M$\;
\
\caption{The local search algorithm of Gotthilf and Lewenstein.}\label{alg2}
\end{algorithm}

\begin{lemma}\label{lemma2}
Let $G$ be a $\{ C_3,C_5\}$-free graph of maximum degree at most $d$ for some $d\geq 3$
such that $\min\{ c_G(e):e\in E(G)\}>\frac{17}{12}d^2$.
If $M$ is an induced matching of $G$ produced by {\sc Local Search} applied to $G$,
then, for every edge $e\in M$,
$pc_G(M,e)\leq \frac{5}{6}d^2+1.$
\end{lemma}
{\it Proof:}
Since $M$ is produced by {\sc Local Search}, it has the following properties.
\begin{enumerate}[(a)]
\item For every edge $e$ of $G$, 
there is some edge $e'$ in $M$ with $e\in C_G(e')$;
because otherwise $e\not\in M$ and $M\cup \{ e\}$ is an induced matching of $G$.
\item If $e$ is in $M$ and $e'$ and $e''$ are two distinct edges in $PC_G(M,e)$,
then $e''\in C_G(e')$;
because otherwise $e',e''\not\in M$ and $(M\setminus \{ e\})\cup \{ e',e''\}$ is an induced matching of $G$.
\end{enumerate}
Let $e=xy\in M$.

Let $X$ be the set of neighbors $u$ of $x$ distinct from $y$ such that $xu\in PC_G(M,e)$.
Let $Y$ be the set of neighbors $u'$ of $y$ distinct from $x$ such that $yu'\in PC_G(M,e)$.
Since $G$ is $\{ C_3,C_5\}$-free, the sets $X$ and $Y$ are disjoint.
Let $X_2$ be the set of vertices $v$ in $V(G)\setminus (\{ x,y\}\cup X\cup Y)$ 
such that there is some vertex $u$ in $X$ with $uv\in PC_G(M,e)$.
Let $Y_2$ be the set of vertices $v'$ in $V(G)\setminus (\{ x,y\}\cup X\cup Y)$ 
such that there is some vertex $u'$ in $Y$ with $u'v'\in PC_G(M,e)$.
Since $G$ is $\{ C_3,C_5\}$-free, the sets $X_2$ and $Y_2$ are disjoint.
Let $X_1$ be the set of vertices in $X$ that have a neighbor in $X_2$.
Let $Y_1$ be the set of vertices in $Y$ that have a neighbor in $Y_2$.
Since $G$ is $\{ C_3,C_5\}$-free, 
the sets $X$, $Y$, $X_2$, and $Y_2$ are independent,
and 
there are no edges between $X_2$ and $Y$
as well as between $Y_2$ and $X$.

By definition, $\{ x,y\}\cup E_G(\{x\},X)\cup E_G(\{y\},Y)$ is the set of edges in $PC_G(M,e)$
that are identical or adjacent with $e$.

Note that is $f\in PC_G(M,e)$ is not identical or adjacent with $e$,  
then $G$ contains an edge $g$ that is adjacent with $e$ and $f$.
If $g\not\in PC_G(M,e)$, then $g\in C_G(e')$ for some edge $e'$ in $M$ distinct from $e$,
which implies the contradiction $f\in C_G(e')$. Hence $g\in PC_G(M,e)$.
This implies that all edges in $PC_G(M,e)\setminus \{ e\}$
are incident with a vertex in $X\cup Y$.

If $uv$ is an edge such that $u\in X_1$ and $v\in X_2$, 
then, by definition, $xu\in PC_G(M,e)$ and $u'v\in PC_G(M,e)$ for some $u'\in X$.
Clearly, $uv\in C_G(e)$.
If $uv\not\in PC_G(M,e)$, then $uv\in C_G(e')$ for some edge $e'$ in $M$ distinct from $e$,
which implies the contradiction that one of the two edges $xu$ and $u'v$ belongs to $C_G(e')$.
Hence, by symmetry, all edges in $E_G(X_1,X_2)\cup E_G(Y_1,Y_2)$ belong to $PC_G(M,e)$.
Similarly, it follows that all edges in $E_G(X,Y)$ belong to $PC_G(M,e)$, and hence
\begin{eqnarray}\label{e1}
PC_G(M,e) & = & \{ x,y\}
\cup E_G(\{x\},X)
\cup E_G(\{y\},Y)\nonumber\\
&&\cup E_G(X_1,X_2)
\cup E_G(Y_1,Y_2)
\cup E_G(X,Y).
\end{eqnarray}
If $u\in X_1$ and $u'\in Y$,
then there is some $v\in X_2$ such that $uv\in PC_G(M,e)$ and, by definition, $yu'\in PC_G(M,e)$.
Since $G$ is $\{ C_3,C_5\}$-free, property (b) implies that $u$ and $u'$ are adjacent,
that is, 
every vertex in $X_1$ is adjacent to every vertex in $Y$,
and, by symmetry,
every vertex in $Y_1$ is adjacent to every vertex in $X$.

If $u_1,u_2\in X_1$ and $v_1,v_2\in X_2$ are four distinct vertices 
such that $u_1v_1$ and $u_2v_2$
are edges of $G$, then, by (\ref{e1}),
$u_1v_1,u_2v_2\in PC_G(M,e)$ and, since $G$ is $\{ C_3,C_5\}$-free, 
property (b) implies that $v_2$ is a neighbor of $u_1$
or $v_1$ is a neighbor of $u_2$,
that is, the bipartite graph between $X_1$ and $X_2$ is $2K_2$-free.
This implies that the sets
$N_G(u)\cap X_2$ for $u$ in $X_1$ are ordered by inclusion.
Hence, 
if $X_1$ is non-empty, 
then $X_1$ contains a vertex $u_x$ that is adjacent to all vertices in $X_2$.
By symmetry, if $Y_1$ is non-empty, 
then the set $Y_1$ contains a vertex $u_y$ that is adjacent to all vertices in $Y_2$.

We consider three cases.

\bigskip

\noindent {\bf Case 1 }{\it $X_1$ and $Y_1$ are both non-empty.}

\bigskip

\noindent Since $u_x$ is adjacent to each vertex in $\{ x\}\cup X_2\cup Y$, 
we obtain 
$|X_2|+|Y|\leq d-1$, 
and, by symmetry,
$|Y_2|+|X|\leq d-1$.
Now (\ref{e1}) implies 
\begin{eqnarray}
pc_G(M,e) & = &
1+m_G(\{x\},X)+m_G(\{y\},Y)+\left( m_G(\{ u_x\},X_2)+m_G(X_1\setminus \{ u_x\},X_2)\right)\nonumber \\
&&+\left( m_G(\{ u_y\},Y_2)+m_G(Y_1\setminus \{ u_y\},Y_2)\right)+m_G(X,Y)\nonumber \\
& = & 1+|X|+|Y|+\left(|X_2|+m_G(X_1\setminus \{ u_x\},X_2)\right)\nonumber\\&&
+\left(|Y_2|+m_G(Y_1\setminus \{ u_y\},Y_2)\right)
+m_G(X,Y)\nonumber \\
& \leq & 2d-1
+m_G(X_1\setminus \{ u_x\},X_2)
+m_G(Y_1\setminus \{ u_y\},Y_2)
+m_G(X,Y)\label{e2}
\end{eqnarray}
Note that $u_x$ is adjacent to all vertices in $X_2\cup Y$ 
and that $x$ is adjacent to all vertices in $\{ y\}\cup (X\setminus \{ u_x\})$.
There are 
$|Y|$ edges between $y$ and $Y$,
$m_G(X,Y)-|Y|$ edges between $X\setminus \{ u_x\}$ and $Y$, and
$m_G(X_1\setminus \{ u_x\},X_2)$ edges between $X_2$ and $X\setminus \{ u_x\}$.
This implies
\begin{eqnarray*}
c_G(xu_x) & \leq & 2d^2-2d+1-m_G(N_G(x)\setminus \{ u_x\},N_G(u_x)\setminus \{ x\})\\
& \leq & 2d^2-2d+1-|Y|-(m_G(X,Y)-|Y|)-m_G(X_1\setminus \{ u_x\},X_2)\\
& = & 2d^2-2d+1-(m_G(X,Y)+m_G(X_1\setminus \{ u_x\},X_2)).
\end{eqnarray*}
Since, by assumption, $c_G(xu_x)>\frac{17}{12}d^2$, 
we obtain
\begin{eqnarray}\label{e3}
m_G(X,Y)+m_G(X_1\setminus \{ u_x\},X_2)&\leq & \frac{7}{12}d^2-2d+1\mbox{, and, by symmetry,}\\
m_G(X,Y)+m_G(Y_1\setminus \{ u_y\},Y_2)&\leq &\frac{7}{12}d^2-2d+1\nonumber.
\end{eqnarray}
By symmetry, we may assume that $|X_1|\geq |Y_1|$.

Since every vertex in $Y_1$ is adjacent to $y$ and to all vertices in $X_1$,
it has at most $d-1-|X_1|$ neighbors in $Y_2$,
which implies 
$m_G(Y_1\setminus \{ u_y\},Y_2)\leq (|Y_1|-1)(d-1-|X_1|)\leq (|Y_1|-1)(d-1-|Y_1|)$.
Note that regardless of the value of $|Y_1|$, 
we have $-|Y_1|^2+d|Y_1|\leq \frac{d^2}{4}$.
Together with (\ref{e2}) and (\ref{e3}) we obtain
\begin{eqnarray*}
pc_G(M,e) 
& \leq & 2d-1
+m_G(X_1\setminus \{ u_x\},X_2)
+m_G(Y_1\setminus \{ u_y\},Y_2)
+m_G(X,Y)\\
& \leq & 2d-1
+\frac{7}{12}d^2-2d+1
+(|Y_1|-1)(d-1-|Y_1|)\\
& = & \frac{7}{12}d^2
-|Y_1|^2+d|Y_1|-d+1\\
& \leq & \frac{7}{12}d^2
+\frac{1}{4}d^2-d+1\\
& = & \frac{5}{6}d^2-d+1.
\end{eqnarray*}
{\bf Case 2 }{\it $X_1$ is non-empty but $Y_1$ is empty.}

\bigskip

\noindent As in Case 1, we have $|X_2|+|Y|\leq d-1$.
Clearly, $|Y|\leq d-1$.
Now (\ref{e1}) implies 
\begin{eqnarray*}
pc_G(M,e) & = &
1
+m_G(\{x\},X)
+m_G(\{y\},Y) +m_G(\{ u_x\},X_2)+m_G(X_1\setminus \{ u_x\},X_2)\\
&&+m_G(X,Y)\nonumber\\
& = & 1+|X|+|Y|+|X_2|
+m_G(X_1\setminus \{ u_x\},X_2)
+m_G(X,Y)\nonumber \\
& \leq & 2d-1
+m_G(X_1\setminus \{ u_x\},X_2)
+m_G(X,Y).
\end{eqnarray*}
Exactly as Case 1, we obtain
\begin{eqnarray*}
m_G(X,Y)+m_G(X_1\setminus \{ u_x\},X_2)&\leq &\frac{7}{12}d^2-2d+1,
\end{eqnarray*}
and hence
\begin{eqnarray*}
pc_G(M,e) 
& \leq & 2d-1
+m_G(X_1\setminus \{ u_x\},X_2)
+m_G(X,Y)\\
& \leq & 2d-1+\frac{7}{12}d^2-2d+1\\
& \leq & \frac{7}{12}d^2.
\end{eqnarray*}
Since $d\geq 3$, it follows that 
$pc_G(M,e)\leq \frac{5}{6}d^2-d+1$.

\bigskip

\noindent {\bf Case 3 }{\it $X_1$ and $Y_1$ are both empty.}

\bigskip

\noindent Note that in this case also both sets $X_2$ and $Y_2$ are empty.
Now (\ref{e1}) implies 
\begin{eqnarray*}
pc_G(M,e) & = &
1
+m_G(\{x\},X)
+m_G(\{y\},Y)
+m_G(X,Y)\\
& \leq & 2d-1+m_G(X,Y).
\end{eqnarray*}
If $m_G(X,Y)\leq \frac{5}{6}d^2-2d+2$, then $pc_G(M,e)\leq \frac{5}{6}d^2+1$.
Hence, we may assume that $m_G(X,Y)>\frac{5}{6}d^2-2d+2$.
Since $|X|\leq d-1$, this implies the existence of a vertex $\tilde{u}_x$ in $X$ 
with at least 
$$\frac{\frac{5}{6}d^2-2d+2}{d-1}
=\frac{5}{6}d-\frac{7}{6}+\frac{5}{6(d-1)}
\geq \frac{5}{6}d-\frac{7}{6}$$
neighbors in $Y$.
Let $\tilde{Y}$ be the set of neighbors of $\tilde{u}_x$ in $Y$.
Since $|Y|\leq d-1$, 
we have $|Y\setminus \tilde{Y}|\leq (d-1)-\left(\frac{5}{6}d-\frac{7}{6}\right)=\frac{1}{6}d+\frac{1}{6}$
and thus
\begin{eqnarray*}
m_G(X\setminus \{ \tilde{u}_x\},\tilde{Y})  & = & 
m_G(X,Y)-m_G(\{\tilde{u}_x\},Y)-m_G(X\setminus \{ \tilde{u}_x\},Y\setminus \tilde{Y})\\
& \geq &
\left(\frac{5}{6}d^2-2d+2\right)-(d-1)-(d-2)\left(\frac{1}{6}d+\frac{1}{6}\right)\\
&=& \frac{2}{3}d^2-\frac{17}{6}d+\frac{10}{3}.
\end{eqnarray*}
Note that $\tilde{u}_x$ is adjacent to all vertices in $\tilde{Y}$ 
and that $x$ is adjacent to all vertices in $\{ y\}\cup (X\setminus \{ \tilde{u}_x\})$.
There are 
$|\tilde{Y}|$ edges between $y$ and $\tilde{Y}$, and
$m_G(X\setminus \{ \tilde{u}_x\},\tilde{Y})$ edges between $X\setminus \{ \tilde{u}_x\}$ and $\tilde{Y}$.
This implies
\begin{eqnarray*}
c_G(x\tilde{u}_x) & \leq & 2d^2-2d+1-m_G(N_G(x)\setminus \{ \tilde{u}_x\},N_G(\tilde{u}_x)\setminus \{ x\})\\
& \leq & 2d^2-2d+1-|\tilde{Y}|-m_G(X\setminus \{ \tilde{u}_x\},\tilde{Y})\\
& \leq & 2d^2-2d+1-\left(\frac{5}{6}d-\frac{7}{6}\right)-\left(\frac{2}{3}d^2-\frac{17}{6}d+\frac{10}{3}\right)\\
& \leq & \frac{4}{3}d^2-\frac{7}{6},
\end{eqnarray*}
contradicting the assumption $c_G(xu_x)>\frac{17}{12}d^2$.

This completes the proof. $\Box$

\begin{lemma}\label{lemma3}
Let $G$ be a $\{ C_3,C_5\}$-free graph of maximum degree at most $d$ for some $d\geq 3$
such that $\min\{ c_G(e):e\in E(G)\}>\frac{17}{12}d^2$.
If $M$ is an induced matching of $G$ produced by {\sc Local Search} applied to $G$,
then
$|M|\geq \frac{m(G)}{\frac{17}{12}d^2-d+1}.$
\end{lemma}
{\it Proof:} We consider the number $p$ of pairs $(e,f)$ where 
$e\in M$ and $f\in C_G(e)$.
Since $c_G(e)\leq 2d^2-2d+1$, we have $p\leq (2d^2-2d+1)|M|$.
In order to obtain a lower bound on $p$,
we observe that the edges $f$ of $G$ for which there is exactly one edge $e$ in $M$
with $f\in C_G(e)$ are exactly those in $\bigcup_{e\in M}PC_G(M,e)$.
Hence, by Lemma \ref{lemma2}, 
$p\geq 2m(G)-\sum_{e\in M}pc_G(M,e)\geq 2m(G)-\left(\frac{5}{6}d^2+1\right)|M|$.
Both estimates for $p$ together yield the desired bound. $\Box$

\bigskip

\begin{algorithm}[H]
{\sc ApproxBip}(d)\\
\KwIn{A $\{ C_3,C_5\}$-free $d$-regular graph $G$.}
\KwOut{An induced matching $M\cup M'$ of $G$.}
\BlankLine
Apply {\sc Greedy$\left(\frac{17}{12}d^2\right)$} to $G$ and denote the output by $(M,G')$\;
Apply {\sc Local Search} to $G'$ and denote the output by $M'$\;
\Return $M\cup M'$\;
\
\caption{An approximation algorithm for $\{ C_3,C_5\}$-free $d$-regular graphs.}\label{alg3}
\end{algorithm}

\begin{lemma}\label{lemma3b}
Let $G$ be a $\{ C_3,C_5\}$-free graph of maximum degree at most $d$ for some $d\geq 3$.
The algorithm {\rm{\sc ApproxBip}(d)} (cf. Algorithm \ref{alg3}) applied to $G$
produces an induced matching $M\cup M'$ 
of $G$ with $|M\cup M'|\geq \frac{12m(G)}{17d^2}=\frac{m(G)}{1.41\bar{6}d^2}$.
\end{lemma}
{\it Proof:} 
By Lemma \ref{lemma1}, the set $M\cup M'$ produced by {\sc ApproxBip}(d) applied to $G$
is an induced matching of $G$.
If $m=m(G)-m(G')$, 
then $|M|\geq \frac{m}{\frac{17}{12}d^2}$.
By Lemma \ref{lemma3}, 
$|M'|\geq \frac{m(G')}{\frac{17}{12}d^2-d+1}\geq \frac{m(G)-m}{\frac{17}{12}d^2}$.
Altogether, we obtain 
$|M\cup M'|\geq \frac{m(G)}{\frac{17}{12}d^2}$.
$\Box$

\bigskip

\noindent We are now in a position to prove Theorem \ref{theorem1}.

\bigskip

\noindent {\it Proof of Theorem \ref{theorem1}:}
For the sake of completeness, 
we give the argument for the upper bound.
If $N$ is some induced matching of the $d$-regular graph $G$,
then, for each edge $e$ in $N$, the two vertices incident with $e$ are incident with exactly $2d-1$ edges of $G$.
Since $N$ is induced, these sets of $2d-1$ edges are disjoint, which implies
$|N|\leq \frac{m(G)}{2d-1}$.
By Lemma \ref{lemma3b}, this implies that the performance ratio 
of {\sc ApproxBip}(d) applied to a $\{ C_3,C_5\}$-free $d$-regular graph
is at most $\frac{17d^2}{12(2d-1)}=\frac{17}{24}d+\frac{17d}{48d-24}$,
which completes the proof. $\Box$

\bigskip

\noindent For the proof of Theorem \ref{theorem2}, we need the following lemma.

\begin{lemma}\label{lemma4}
If $G$ is a non-empty $k$-degenerate graph of maximum degree at most $d$ with $k<d$,
then $G$ has an edge $e$ with 
$c_G(e)\leq (3k-1)d-k(k+1)+1.$
\end{lemma}
{\it Proof:} Let $X$ denote the set of vertices of $G$ of degree at most $k$ 
and let $Y=V(G)\setminus X$.
If two vertices $u$ and $v$ in $X$ are adjacent, 
then $c_G(e)\leq 1+2(k-1)d$.
If $X$ is an independent set, 
then, since the graph $G[Y]$ has a vertex $v$ of degree at most $k$ in $G[Y]$ but degree more than $k$ in $G$,
some vertex $u$ in $X$ is adjacent to this vertex $v$ in $Y$.
Each of the at most $k$ neighbors of $u$ has degree at most $d$.
Since $v$ has degree at most $k$ in $G[Y]$, 
it has $r$ neighbors in $Y$ for some $r\leq k$,
which have degree at most $d$.
The remaining at most $d-r$ neighbors of $v$ are all in $X$
and hence have degree at most $k$.
Altogether this implies that 
\begin{eqnarray*}
c_G(uv) & \leq & 1+(k-1)d+rd+(d-1-r)k\\
&= & 1+2kd+r(d-k)-d-k\\
&\leq & 1+2kd+k(d-k)-d-k\\
&= & 1+(3k-1)d-k(k+1).
\end{eqnarray*}
Since $2(k-1)d\leq (3k-1)d-k(k+1)$, the proof is complete.
$\Box$

\bigskip

\noindent {\it Proof of Theorem \ref{theorem2}:} 
By definition, every subgraph of $G$ is $k$-degenerate and of maximum degree at most $d$.
Therefore, by Lemma \ref{lemma4},
the algorithm {\sc \sc Greedy($(3k-1)d-k(k+1)+1$)} applied to the graph $G$ 
outputs a pair $(M,G')$ where $M$ is an induced matching of $G$
and the graph $G'$ has no edges.
This implies that $|M|\geq \frac{m(G)}{(3k-1)d-k(k+1)+1}$,
which completes the proof. $\Box$


\begin{thebibliography}{}

\bibitem{brmo}
A. Brandst\"{a}dt and R. Mosca, On distance-3 matchings and induced matchings, Discrete Appl. Math. 159 (2011) 509-520.
\bibitem{ca}
K. Cameron, Induced matchings, Discrete Appl. Math. 24 (1989) 97-102.
\bibitem{ca2}
K. Cameron, Induced matchings in intersection graphs, Discrete Math. 278 (2004) 1-9.
\bibitem{cn}
G.J. Chang and N. Narayanan, Strong chromatic index of 2-degenerate graphs, J. Graph Theory 73 (2013) 119-126.
\bibitem{dadelo}
K.K. Dabrowski, M. Demange, and V.V. Lozin, New results on maximum induced matchings in bipartite graphs and beyond, 
Theor. Comput. Sci. 478 (2013) 33-40.
\bibitem{dgs}
M. D\c{e}bski, J. Grytczuk and M. \'{S}leszy\'{n}ska-Nowak, 
Strong chromatic index of sparse graphs, 
arXiv:1301.1992 (9 January 2013).

\bibitem{dumazi}
W. Duckworth, D.F. Manlove, and M. Zito,
On the approximability of the maximum induced matching problem,
J. Discrete Algorithms 3 (2005) 79-91.

\bibitem{fj}
J.L. Fouquet and J.L. Jolivet, 
Strong edge-colouring of graphs and applications to multi-$k$-gons,
Ars Combin. 16A (1983) 141-150.

\bibitem{gole}
Z. Gotthilf and M. Lewenstein,
Tighter Approximations for Maximum Induced Matchings in Regular Graphs,
Lecture Notes in Computer Science 3879 (2006) 270-281.

\bibitem{jorasa}
F. Joos, D. Rautenbach, and T. Sasse,
Induced Matchings in Subcubic Graphs,
to appear in SIAM J. Discrete Math.

\bibitem{lopl}
L. Lov\'{a}sz and M.D. Plummer,
Matching Theory, vol. 29, Annals of Discrete Mathematics, North-Holland, Amsterdam, 1986.

\bibitem{lo}
V.V. Lozin,
On maximum induced matchings in bipartite graphs,
Inf. Process. Lett. 81 (2002) 7-11.

\bibitem{ly}
R. Luo and G. Yu, 
A note on strong edge-colorings of $2$-degenerate graphs,
arXiv:1212.6092 (25 December 2012).

\bibitem{more} 
M. Molloy and B. Reed,
A bound on the strong chromatic index of a graph, 
J. Combin. Theory Ser. B 69 (1997) 103-109.

\bibitem{stva}
L.J. Stockmeyer and V.V. Vazirani,
NP-completeness of some generalizations of the maximum matching problem,
Inf. Process. Lett. 15 (1982) 14-19.

\bibitem{y}
G. Yu,
Strong edge-colorings for $k$-degenerate graphs,
arXiv:1212.6093 (31 March 2013).

\bibitem{w}
T. Wang,
Strong chromatic index of k-degenerate graphs,
arXiv:1304.0285 (19 July 2013).

\bibitem{zi}
M. Zito,
Induced Matchings in Regular Graphs and Trees, 
Lecture Notes in Computer Science 1665 (1999) 89-100.

\end{thebibliography}
\end{document}